\documentclass[leqno, 11pt]{amsart}
\usepackage{amsmath}
\usepackage{amsfonts}
\usepackage{amssymb}
\usepackage{graphicx}%
\newtheorem{theorem}{Theorem}

\newtheorem{corollary}[theorem]{Corollary}

\newtheorem{lemma}[theorem]{Lemma}

\newtheorem{problem}[theorem]{Problem}

\newtheorem{remark}[theorem]{Remark}

\begin{document}

\title[Characterizations of model manifolds]
{Characterizations of model manifolds by means of certain differential systems}

\begin{abstract}
We prove metric rigidity for complete manifolds supporting solutions of
certain second order differential systems, thus extending classical works on a
characterization of space-forms. In the route, we also discover
new characterizations of space-forms. We next generalize results concerning metric
rigidity via equations involving vector fields.

\end{abstract}

\date{March 05, 2009}
\author{S. Pigola}
\address{Dipartimento di Fisica e Matematica\\Via Valleggio 11, 22100 Como, Italy}
\email{stefano.pigola@uninsubria.it}
\author{M. Rimoldi}
\address{Dipartimento di Matematica\\Via Saldini 50, 20133 Milano, Italy}
\email{michele.rimoldi@unimi.it}

\maketitle

\section{Introduction\label{section-introduction}}

Having fixed a smooth, even function $G:\mathbb{R}\rightarrow\mathbb{R}$, we
let $M_{-G}^{m}$ denote the $m$-dimensional (not necessarily complete) model
manifold with radial sectional curvature $-G\left(  r\right)  $. More
precisely, we set%
\[
M_{-G}^{m}=\left(  [0,r_{-G})\times\mathbf{S}^{m-1},dr^{2}+g\left(  r\right)
^{2}d\theta^{2}\right)  ,
\]
where $g:\mathbb{R}\rightarrow\mathbb{R}$ is the unique solution of the
problem%
\[
\left\{
\begin{array}
[c]{l}%
g^{\prime\prime}=Gg\\
g\left(  0\right)  =0\\
g^{\prime}\left(  0\right)  =1,
\end{array}
\right.
\]
and $r_{-G}\in(0,+\infty]$ is the first zero of $g\left(  r\right)  $ on
$(0,+\infty)$. Obviously, in case $g\left(  r\right)  >0$ for every $r>0$, we
are using the convention $r_{-G}=+\infty$. In this case, the model is
geodesically complete.

Examples of models comes from the standard space-forms.

\begin{enumerate}
\item[(a)] Let $G\left(  r\right)  \equiv-k<0$. Then $g\left(  r\right)
=k^{-1/2}\sin\left(  k^{1/2}r\right)  $, $r_{k}=\pi/k^{1/2}$ and $M_{k}^{m}$
is isometric to the standard sphere of constant curvature $k$ punctured at one
point. Equivalently, $M_{k}^{m}$ is isometric to the geodesic ball
$B_{\pi/\sqrt{k}}\left(  o\right)  $ in the standard sphere of constant
curvature $k$.

\item[(b)] Let $G\left(  r\right)  \equiv k>0$. Then $g\left(  r\right)
=k^{-1/2}\sinh\left(  k^{1/2}r\right)  $ and $M_{-k}^{m}$ is isometric to the
standard hyperbolic space of constant curvature $-k.$

\item[(c)] Let $G\left(  r\right)  \equiv0$. Then $g\left(  r\right)  =r$ and
$M_{0}^{m}$ is isometric to the standard Euclidean space.
\end{enumerate}

Characterizations of space-forms as complete manifolds supporting solutions of
second order differential systems of the form%
\[
\mathrm{Hess}\left(  u\right)  \left(  x\right)  =\left(  au\left(  x\right)
+b\right)  \left\langle ,\right\rangle _{x},
\]
have been classically investigated by M. Obata, \cite{Obata-sphere}, Y.
Tashiro, \cite{Tashiro-vector}, and M. Kanai, \cite{Kanai-structure}. The
following theorem encloses in a single statement their results.

\begin{theorem}
\label{th_spaceforms}Let $\left(  M,\left\langle ,\right\rangle \right)  $ be
a complete, connected Riemannian manifold of dimension $\dim M=m$. Then:

\begin{enumerate}
\item[(a)] A necessary and sufficient condition for $M$ to be isometric to the
sphere of constant curvature $k>0$ is that $M$ supports a smooth, non trivial
solution $u:M\rightarrow\mathbb{R}$ of the differential system%
\begin{equation}
\mathrm{Hess}\left(  u\right)  \left(  x\right)  =-ku\left(  x\right)
\left\langle ,\right\rangle . \label{spaceforms-sphere}%
\end{equation}

\item[(b)] A necessary and sufficient condition for $M$ to be isometric to the
hyperbolic space of constant curvature $-k<0$ is that $M$ supports a smooth,
non trivial solution $u:M\rightarrow\mathbb{R}$ of the differential system%
\begin{equation}
\mathrm{Hess}\left(  u\right)  \left(  x\right)  =ku\left(  x\right)
\left\langle ,\right\rangle , \label{spaceforms-hyperbolic}%
\end{equation}
with precisely one critical point.

\item[(c)] A necessary and sufficient condition for $M$ to be isometric to the
Euclidean space is that $M$ supports a smooth, non trivial solution
$u:M\rightarrow\mathbb{R}$ of the differential system%
\begin{equation}
\mathrm{Hess}\left(  u\right)  \left(  x\right)  =h\left\langle ,\right\rangle
, \label{spaceforms-euclidean}%
\end{equation}
for some constant $h\neq0$.
\end{enumerate}
\end{theorem}

Recently, E. Garcia-Rio, D. Kupeli and B. Unal, \cite{GKU-spheres}, have been
able to extend the metric rigidity established in Theorem \ref{th_spaceforms}
to complete manifolds supporting vector field solutions $Z$ of differential
systems of the form%
\[
\left(  DDZ\right)  \left(  X,Y\right)  =k\left\langle Z,X\right\rangle Y,
\]
for some constant $k\neq0$ and for every vector fields $X,Y$. Here, the symbol
$D$ stands for covariant differentiation so that%
\[
\left(  DDZ\right)  \left(  X,Y\right)  =D_{X}D_{Y}Z-D_{D_{X}Y}Z.
\]
Note that, in case $Z=\nabla u$ is a gradient vector field, the above equation
becomes%
\[
D\mathrm{Hess}\left(  u\right)  =k\text{ }du\otimes\left\langle ,\right\rangle
,
\]
which is a third order system in the unknown function $u$. The following
rigidity theorem summarizes the main results of \cite{GKU-spheres}.

\begin{theorem}
\label{th_spaceforms-vector}Let $\left(  M,\left\langle ,\right\rangle
\right)  $ be a complete, connected Riemannian manifold of dimension $\dim
M=m$. Then:

\begin{enumerate}
\item[(a)] A necessary and sufficient condition for $M$ to be isometric to the
sphere of constant curvature $k>0$ is that $M$ supports a smooth, non trivial
solution $Z$ of the differential system%
\[
\left(  DDZ\right)  \left(  X,Y\right)  =-k\left\langle Z,X\right\rangle
Y,\text{\qquad}\forall X,Y.
\]

\item[(b)] A necessary and sufficient condition for $M$ to be isometric to the
hyperbolic space of constant curvature $-k<0$ is that $M$ supports a smooth,
non trivial solution $Z$ of the differential system%
\[
\left(  DDZ\right)  \left(  X,Y\right)  =k\left\langle Z,X\right\rangle
Y,\text{\qquad}\forall X,Y
\]
satisfying $Z_{o}=0$, for some $o\in M$.
\end{enumerate}
\end{theorem}

Since space-forms are very special cases of model manifolds, a natural
question is whether a general model manifold $M_{-G}^{m}$ can be characterized
in the same perspective of Theorem \ref{th_spaceforms} and Theorem
\ref{th_spaceforms-vector}. This note aims to answer the question in the affirmative.
During our investigation, we will also give new characterizations of space-forms.

\section{Second order systems}

Quite naturally, one expects that a characterization of the model $M_{-G}^{m}%
$, in the spirit of Theorem \ref{th_spaceforms} above, must involve more
general differential systems of the form%
\begin{equation}
\mathrm{Hess}\left(  u\right)  \left(  x\right)  =H\left(  r\left(  x\right)
\right)  u\left(  x\right)  \left\langle ,\right\rangle , \label{radial}%
\end{equation}
where $r\left(  x\right)  $ denotes the geodesic distance from a fixed origin
$o$. First of all, we need to find the right form of the radial coefficient
$H$. Let $u\left(  x\right)  =\alpha\left(  r\left(  x\right)  \right)  $ be a
radial solution of (\ref{radial}). We suppose to have normalized $u$ in such a
way that $u\left(  0\right)  =1$ and we require $u$ to have a critical point
at $0$. Then, recalling that%
\begin{equation}
\mathrm{Hess}\left(  r\right)  =\frac{g^{\prime}}{g}\left\{  \left\langle
,\right\rangle -dr\otimes dr\right\}  =gg^{\prime}d\theta^{2}, \label{radial2}%
\end{equation}
we have%
\[
\mathrm{Hess}\left(  u\right)  =\alpha^{\prime\prime}dr\otimes dr+\alpha
^{\prime}gg^{\prime}d\theta^{2}%
\]
On the other hand%
\[
\mathrm{Hess}\left(  u\right)  =H\alpha\left\langle ,\right\rangle =H\alpha
dr\otimes dr+H\alpha g^{2}d\theta^{2}.
\]
Comparing these two equations gives the ordinary differential system%
\[
\left\{
\begin{array}
[c]{l}%
\alpha^{\prime\prime}=H\alpha\\
\alpha^{\prime}gg^{\prime}=H\alpha g^{2},
\end{array}
\right.
\]
that is,%
\[
\left\{
\begin{array}
[c]{l}%
\alpha^{\prime\prime}=\alpha^{\prime}g^{\prime}/g\\
H=\alpha^{\prime}g^{\prime}/\alpha g,
\end{array}
\right.
\]
where, we recall, $\alpha\left(  0\right)  =1$, $\alpha^{\prime}\left(
0\right)  =0$. Integrating the first equation gives%
\begin{equation}
\alpha\left(  r\right)  =A\int_{0}^{r}g\left(  s\right)  ds+1, \label{radial3}%
\end{equation}
with $A\neq0$ any constant. Inserting this expression into the second equation
we finally deduce%
\[
H\left(  r\right)  =\frac{Ag^{\prime}\left(  r\right)  }{A\int_{0}^{r}g\left(
s\right)  ds+1}.
\]
In order that $H$ is defined on all of $[0,r_{-G})$ we need to impose that%
\[
\inf\left\{  t>0:A\int_{0}^{t}g\left(  s\right)  ds+1\leq0\right\}  \geq
r_{-G}.
\]
We have thus obtained the following

\begin{lemma}
A necessary and sufficient condition for equation (\ref{radial}) on
$M_{-G}^{m}$ to possess a radial solution $u$ is that%
\[
H\left(  r\right)  =\frac{Ag^{\prime}\left(  r\right)  }{A\int_{0}^{r}g\left(
s\right)  ds+1}.
\]
for any constant $A\neq0$ such that%
\[
\inf\left\{  t>0:A\int_{0}^{t}g\left(  s\right)  ds+1\leq0\right\}  \geq
r_{-G}.
\]

\end{lemma}

Note that, in particular,

\begin{enumerate}
\item[(a)] On the punctured standard sphere $M_{1}^{m}=\mathbf{S}%
^{m}\backslash\left\{  \text{point}\right\}  =B_{\pi}\left(  0\right)  $, for
every $A\in\mathbb{R}\backslash\left\{  0\right\}  $ such that $A>-1/2,$
$A=-1$, there is a smooth function $u_{A}$ with exactly one critical point at
$0$ and satisfying the equation%
\begin{equation}
\mathrm{Hess}\left(  u_{A}\right)  \left(  x\right)  =\frac{A\cos r\left(
x\right)  }{-A\cos r\left(  x\right)  +1+A}u_{A}\left(  x\right)  \left\langle
,\right\rangle . \label{spaceforms2}%
\end{equation}
As a matter of fact, the function $u\left(  x\right)  =-A\cos r\left(
x\right)  +1+A$ is well defined and solves the equation on all of
$\mathbf{S}^{m}$. Note finally that, in the special case $A=-1,$
(\ref{spaceforms2}) reduces to (\ref{spaceforms-sphere}).

\item[(b)] On the standard hyperbolic model $M_{-1}^{m}=\mathbb{H}_{-1}^{m}$,
for every $A>0$, there exists a smooth function $u_{A}$ with exactly one
critical point at $0$ and satisfying the equation%
\begin{equation}
\mathrm{Hess}\left(  u_{A}\right)  \left(  x\right)  =\frac{A\cosh r\left(
x\right)  }{A\cosh r\left(  x\right)  +1-A}u_{A}\left(  x\right)  \left\langle
,\right\rangle . \label{spaceforms1}%
\end{equation}
In the special case $A=1$, (\ref{spaceforms1}) reduces to
(\ref{spaceforms-hyperbolic}).

\item[(c)] On the standard Euclidean space $M_{0}^{m}=\mathbb{R}^{m},$ for
every $A>0$, there exists a function $u_{A}$ with exactly one critical point
at $0$ and satisfying the equation%
\begin{equation}
\mathrm{Hess}\left(  u_{A}\right)  \left(  x\right)  =\frac{2A}{Ar\left(
x\right)  ^{2}+2}u_{A}\left(  x\right)  \left\langle ,\right\rangle .
\label{spaceforms3}%
\end{equation}

\end{enumerate}

We shall prove the following result. Recall that a twisted sphere of dimension
$n$ is a differentiable manifold $N$, homeomorphic to the standard sphere
$\mathbf{S}^{n}$, which is obtained by gluing two $n$-dimensional closed, unit
disks $D^{n}\subset\mathbb{R}^{n}$ via a boundary diffeomorphism.

\begin{theorem}
\label{th_global}Let $\left(  M,\left\langle ,\right\rangle \right)  $ be a
complete Riemannian \ manifold of dimension $m$, and let $o\in M$ be a
reference origin. Then, a necessary and sufficient condition for the existence
of an isometric imbedding $\Phi:M_{-G}^{m}\rightarrow M$ is that there exists
a smooth solution $u:B_{r_{-G}}\left(  o\right)  \rightarrow\mathbb{R}$ of the
problem%
\begin{equation}
\left\{
\begin{array}
[c]{l}%
\mathrm{Hess}\left(  u\right)  \left(  x\right)  =H\left(  r\left(  x\right)
\right)  u\left(  x\right)  \left\langle ,\right\rangle \\
u\left(  o\right)  =1\\
\left\vert \nabla u\right\vert \left(  o\right)  =0,
\end{array}
\right.  \label{local1}%
\end{equation}
where $r\left(  x\right)  =\mathrm{dist}_{\left(  M,\left\langle
,\right\rangle \right)  }\left(  x,o\right)  $, $H:[0,R^{\ast}]\rightarrow
\mathbb{R}$ is the smooth function%
\begin{equation}
H\left(  t\right)  =\frac{Ag^{\prime}\left(  t\right)  }{A\int_{0}^{t}g\left(
s\right)  ds+1}, \label{local2}%
\end{equation}
for some real number $A\neq0$, and%
\[
R^{\ast}=\sup\left\{  T>0:H\left(  t\right)  \text{ well defined on
}[0,T]\right\}  >r_{-G}.
\]
Furthermore, if $u$ is a solution of (\ref{local1}) on all of $M$, then the
following holds:

\begin{enumerate}
\item[(a)] In case $r_{-G}=+\infty$, then $M$ is isometric to the model
$M_{-G}^{m}$.

\item[(b)] In case $r_{-G}<+\infty$ and $H\left(  r_{-G}\right)  \neq0$, then
$cut\left(  o\right)  =\left\{  O\right\}  $ for some $O\in M$, and
$\Phi\left(  M_{-G}^{m}\right)  =M\backslash\left\{  O\right\}  $.
Furthermore, $M$ is diffeomorphically a twisted sphere.
\end{enumerate}
\end{theorem}

As a direct consequence of Theorem \ref{th_global} we point out the following
result that generalizes, in some directions, Theorem \ref{th_spaceforms} above.

\begin{corollary}
\label{cor_spaceforms}Let $\left(  M,\left\langle ,\right\rangle \right)  $ be
a complete Riemannian manifold, $o\in M$ a reference origin and $r\left(
x\right)  =\mathrm{dist}_{\left(  M,\left\langle ,\right\rangle \right)
}\left(  x,o\right)  $. Then:

\begin{enumerate}
\item[(a)] $M$ is isometric to the standard sphere $\mathbf{S}^{m}$ if and
only if $M$ supports a real valued function $u\not \equiv 0$ with a critical
point at $o$ and satisfying the differential system (\ref{spaceforms2}), for
some $A\neq0$ such that either $A>-1/2$ or $A=-1$.

\item[(b)] $M$ is isometric to the standard hyperbolic space if and only if
$M$ supports a real valued function $u\not \equiv 0$ with a critical point at
$o$ and satisfying the differential system (\ref{spaceforms1}) for some $A>0.$

\item[(c)] $M$ is isometric to the standard Euclidean space if and only if $M$
supports a real valued function $u\not \equiv 0$ with a critical point at $o$
and satisfying the differential system (\ref{spaceforms3}) for some $A>0.$
\end{enumerate}
\end{corollary}

Before proving Theorem \ref{th_global} we make some observations on case (a)
of the previous Corollary.

\begin{enumerate}
\item[(i)] First of all, to deduce that $M$ is a standard sphere one simply
observes that, as established in (b) of Theorem \ref{th_global}, $M$ is simply
connected and $M\backslash\left\{  O\right\}  $ is isometric to a standard
punctured sphere. Therefore, by continuity, $M$ itself has positive constant
curvature and we can apply the Hopf classification theorem. Alternatively, we
can recall that a necessary and sufficient condition for the model metric
$dr\otimes dr+g\left(  r\right)  ^{2}d\theta^{2}$ of $M_{-G}^{m}$ to smoothly
extend on all of $[0,r_{-G}]\times\mathbf{S}^{m-1}$ is that $g^{\left(
2k\right)  }\left(  r_{-G}\right)  =0$ and $g^{\prime}\left(  r_{-G}\right)
=-1$; see \cite{Petersen-riemannian}. In the present situation we have
$g\left(  r\right)  =\sin\left(  r\right)  $ and therefore we deduce that the
isometry $\Phi$ extends to cover the removed point $O$.

\item[(ii)] Comparing with case (a) of Theorem \ref{th_spaceforms} we see
that, on the one hand, we enlarge the class of differential systems
characterizing the sphere but, on the other hand, we make the additional
assumption that $u$ has a critical point at $o$. As first noted by Obata, the
existence of a critical point is automatically guaranteed if $H\left(
r\right)  \equiv-k<0.$ To see this, one can argue as follows. By
contradiction, suppose $u$ has no critical point at all. Then, the vector
field $X=\nabla u/\left\vert \nabla u\right\vert $ is defined on all of $M$.
Using the differential system $\mathrm{Hess}\left(  u\right)  =-ku\left\langle
,\right\rangle $ it is readily seen that the integral curves $\gamma\left(
t\right)  :\mathbb{R}\rightarrow M$ of $X$ are unit speed, but not necessarily
minimizing, geodesics. Indeed%
\begin{align*}
D_{\overset{\cdot}{\gamma}}\overset{\cdot}{\gamma}  &  =D_{\overset{\cdot
}{\gamma}}X_{\gamma}\\
&  =\left\vert \nabla u\right\vert ^{-1}Hess\left(  u\right)  \left(
\overset{\cdot}{\gamma},\cdot\right)  ^{\#}-\left\vert \nabla u\right\vert
^{-1}Hess\left(  u\right)  \left(  \overset{\cdot}{\gamma},X\right)  X\\
&  =-ku\left\vert \nabla u\right\vert ^{-1}X+ku\left\vert \nabla u\right\vert
^{-1}X\\
&  =0.
\end{align*}
Note that the same argument works if $u$ solves the more general equation
$\mathrm{Hess}\left(  u\right)  =f\left\langle ,\right\rangle $, for any
real-valued function $f$. Now consider $y\left(  t\right)  =u\circ
\gamma\left(  t\right)  .$ Then, $y$ satisfies the oscillatory o.d.e.%
\[
y^{\prime\prime}=-ky.
\]
Let $t_{0}>0$ be a critical point of $y$. Since%
\begin{align*}
0  &  =\frac{dy}{dt}\left(  t_{0}\right) \\
&  =\left\langle \nabla u\left(  \gamma\left(  t_{0}\right)  \right)
,\overset{\cdot}{\gamma}\left(  t_{0}\right)  \right\rangle \\
&  =\left\langle \nabla u\left(  \gamma\left(  t_{0}\right)  \right)
,\frac{\nabla u}{\left\vert \nabla u\right\vert }\left(  \gamma\left(
t_{0}\right)  \right)  \right\rangle \\
&  =\left\vert \nabla u\right\vert \left(  \gamma\left(  t_{0}\right)
\right)  ,
\end{align*}
we have that $\gamma\left(  t_{0}\right)  $ is a critical point of $u$.
Contradiction. Thus, $u$ has a critical point $p$ and we can always take $p=o$
as the reference origin in our Theorem \ref{th_global}.

In case the coefficient $H$ in the differential equation depends on the
distance function $r\left(  x\right)  $, if we try to adapt the previous
argument to the present situation, we encounter two obvious difficulties.

\begin{enumerate}
\item As observed above, an integral curve $\gamma\left(  t\right)
:\mathbb{R}\rightarrow M$ of the vector field $X$ is a geodesic but it can be
non-minimizing. Therefore, for large values of $\left\vert t\right\vert $,
$H\left(  r\left(  \gamma\left(  t\right)  \right)  \right)  \neq H\left(
t\right)  .$ It follows that the reduction procedure of the P.D.E. to an
o.d.e., via composition with $\gamma$, cannot be carried over for large values
of $\left\vert t\right\vert $.

\item Even if we were able to prove that $u$ has a critical point at some
$p\in M$, since the coefficient $H$ depends on the distance from the reference
origin $o$, we could not take $p=o$.
\end{enumerate}
\end{enumerate}

The rest of the section is entirely devoted to a proof of Theorem
\ref{th_global}. The \textquotedblleft necessity\textquotedblright\ part has
been already discussed above. Therefore we may concentrate on the
\textquotedblleft sufficiency\textquotedblright\ part.

The following density result due to R. Bishop, \cite{Bishop-cutloci}, will
play a key role in our argument. For a nice and simplified proof, see F.
Wolter, \cite{Wolter-cutloci}. Following Bishop, recall that, given a complete
manifold $\left(  M,\left\langle ,\right\rangle \right)  $ and a reference
point $o\in M$, then $p\in cut\left(  o\right)  $ is an\textit{ ordinary cut
point }if there are at least two distinct minimizing geodesics from $o$ to
$p$. Using the infinitesimal Euclidean law of cosines, it is not difficult to
show that at an ordinary cut point $p$ the distance function $r\left(
x\right)  =\mathrm{dist}_{\left(  M,\left\langle ,\right\rangle \right)
}\left(  x,o\right)  $ is not differentiable, \cite{Wolter-cutloci}.

\begin{theorem}
\label{th_bishop}Let $\left(  M,\left\langle ,\right\rangle \right)  $ be a
complete Riemannian manifold and let $o\in M$ be a reference point. Then the
ordinary cut-points of $o$ are dense in $cut\left(  o\right)  $. In
particular, if the distance function $r\left(  x\right)  $ from $o$ is
differentiable on the (punctured) open ball $B_{R}\left(  o\right)
\backslash\left\{  o\right\}  $ then $B_{R}\left(  o\right)  \cap cut\left(
o\right)  =\emptyset$.
\end{theorem}

We now come into the

\begin{proof}
[Proof (of Theorem \ref{th_global})]To simplify the exposition we will proceed
by steps. The strategy of the proof essentially follows Obata's original
path.\medskip

\textbf{Step 1.} First of all, we note that the function $u:B_{r_{-G}}\left(
o\right)  \rightarrow\mathbb{R}$ must be radial and, more precisely,%
\[
u\left(  x\right)  =\alpha\left(  r\left(  x\right)  \right)  ,
\]
where%
\[
\alpha\left(  t\right)  =A\int_{0}^{t}g\left(  s\right)  ds+1.
\]
Indeed, fix $x$ and choose a unit speed, minimizing geodesic $\gamma
:[0,r\left(  x\right)  ]\rightarrow B_{r_{-G}}\left(  o\right)  $ from $o$ to
$x$. Then, composing with $\gamma$ we deduce that, $y\left(  t\right)
=u\circ\gamma\left(  t\right)  $ is the solution of the Cauchy problem%
\[
\left\{
\begin{array}
[c]{l}%
y^{\prime\prime}\left(  t\right)  =\frac{Ag^{\prime}\left(  t\right)  }%
{A\int_{0}^{t}g\left(  s\right)  ds+1}y\left(  t\right) \\
y\left(  0\right)  =1\\
y^{\prime}\left(  0\right)  =\left\langle \nabla u\left(  o\right)
,\overset{\cdot}{\gamma}\left(  0\right)  \right\rangle =0.
\end{array}
\right.
\]
It follows that%
\[
y\left(  t\right)  =A\int_{0}^{t}g\left(  s\right)  ds+1,
\]
and, taking $t=r\left(  x\right)  $, we get%
\[
u\left(  x\right)  =y\left(  r\left(  x\right)  \right)  =A\int_{0}^{r\left(
x\right)  }g\left(  s\right)  ds+1.
\]
\medskip

\textbf{Step 2.} The open ball $B_{r_{-G}}\left(  o\right)  $ is inside the
cut-locus of $o$. Indeed, recall that $u\left(  x\right)  =\alpha\left(
r\left(  x\right)  \right)  $ and note that $\alpha$ is a diffeomorphism on
$(0,r_{-G})$ because $\alpha^{\prime}\left(  t\right)  =Ag\left(  t\right)
\neq0$ on that interval. Therefore, $r\left(  x\right)  =\alpha^{-1}\circ
u\left(  x\right)  $ is smooth on $B_{r_{-G}}\left(  o\right)  \backslash
\left\{  o\right\}  $ as a composition of smooth functions. By Theorem
\ref{th_bishop}, it follows that $B_{r_{-G}}\left(  o\right)  \cap cut\left(
o\right)  =\emptyset$.\medskip

\textbf{Step 3.} According to Step 2, we can introduce geodesic polar
coordinates on $B_{r_{-G}}\left(  o\right)  $. We claim that the corresponding
map%
\[
\Phi\left(  r,\theta\right)  =\exp_{o}\left(  r\theta\right)  :M_{-G}%
^{m}\approx\mathbf{B}_{r_{-G}}^{m}\left(  0\right)  \subseteq T_{o}%
M\rightarrow B_{r_{-G}}\left(  o\right)  \subseteq M
\]
is a Riemannian isometry. To see this, let $v$ be the function%
\[
v\left(  x\right)  =\frac{u\left(  x\right)  -1}{A}=\int_{0}^{r\left(
x\right)  }g\left(  s\right)  ds
\]
on $B_{r_{-G}}\left(  o\right)  $ and note that%
\begin{equation}
\left\{
\begin{array}
[c]{l}%
\mathrm{Hess}\left(  v\right)  =A^{-1}Hu\left\langle ,\right\rangle \\
v\left(  o\right)  =0\\
\left\vert \nabla v\right\vert \left(  o\right)  =0.
\end{array}
\right.  \label{local3}%
\end{equation}
Furthermore,%
\begin{equation}
\nabla r=\frac{\nabla v}{\left\vert \nabla v\right\vert }. \label{local4}%
\end{equation}
Using geodesic polar coordinates $\left(  r,\theta\right)  \in\left(
0,r_{-G}\right)  \times\mathbf{S}^{m-1}\approx\mathbf{B}_{r_{-G}}^{m}\left(
0\right)  \backslash\left\{  0\right\}  \subseteq T_{o}M$, keeping a local
orthonormal frame $\left\{  \theta^{\alpha}\right\}  $ on $\mathbf{S}%
^{m-1}\subset T_{o}M$, and recalling Gauss lemma, we now express%
\[
\exp_{o}^{\ast}\left\langle ,\right\rangle =dr\otimes dr+\sigma_{\alpha\beta
}\left(  r,\theta\right)  \theta^{\alpha}\otimes\theta^{\beta},
\]
where $d\theta^{2}=\sum\theta^{\alpha}\otimes\theta^{\alpha}$ denotes the
standard metric on $\mathbf{S}^{m-1}$ and the coefficient matrix $\left(
\sigma_{\alpha\beta}\right)  $ satisfies the asymptotic condition%
\begin{equation}
\sigma_{\alpha\beta}\left(  r,\theta\right)  =r^{2}\delta_{\alpha\beta
}+o\left(  r^{2}\right)  ,\text{ as }r\rightarrow0. \label{local4a}%
\end{equation}
By the fundamental equations of Riemannian geometry, we know that, within the
cut locus of $o$,
\[
L_{\nabla r}\left\langle ,\right\rangle =2\mathrm{Hess}\left(  r\right)  ,
\]
where, furthermore, $\nabla r=\partial_{r}$ the radial vector field.
Therefore, on $B_{r_{-G}}\left(  o\right)  $, we have%
\begin{equation}
\partial_{r}\sigma_{\alpha\beta}\left(  r,\theta\right)  =2\mathrm{Hess}%
\left(  r\right)  _{\alpha\beta}. \label{local5}%
\end{equation}
But, according to (\ref{local3}) and (\ref{local4}), we have, for every
$X,Y\in\left(  \nabla r\right)  ^{\bot}$,%
\begin{align*}
\mathrm{Hess}\left(  r\right)  \left(  X,Y\right)   &  =\left\langle
D_{X}\frac{\nabla v}{\left\vert \nabla v\right\vert },Y\right\rangle \\
&  =\frac{1}{\left\vert \nabla v\right\vert }\mathrm{Hess}\left(  v\right)
\left(  X,Y\right) \\
&  =\frac{1}{\left\vert \nabla v\right\vert }A^{-1}Hu\left\langle
X,Y\right\rangle \\
&  =\frac{g^{\prime}}{g}\left\langle X,Y\right\rangle .
\end{align*}
Using this information into (\ref{local5}) and recalling (\ref{local4a}) we
deduce that%
\begin{equation}
\left\{
\begin{array}
[c]{l}%
\partial_{r}\sigma_{\alpha\beta}\left(  r,\theta\right)  =2\dfrac{g^{\prime}%
}{g}\left(  r\right)  \sigma_{\alpha\beta}\left(  r,\theta\right)  \bigskip\\
\sigma_{\alpha\beta}\left(  r,\theta\right)  =r^{2}\delta_{\alpha\beta
}+o\left(  r^{2}\right)  ,\text{ as }r\rightarrow0,
\end{array}
\right.  \label{local6}%
\end{equation}
which integrated gives%
\[
\sigma_{\alpha\beta}\left(  r,\theta\right)  =g\left(  r\right)  ^{2}%
\delta_{\alpha\beta}.
\]
We have thus shown that%
\[
\exp_{o}^{\ast}\left\langle ,\right\rangle =dr\otimes dr+g\left(  r\right)
^{2}d\theta^{2},
\]
proving that $\exp_{o}:M_{-G}^{m}\backslash\left\{  0\right\}  \rightarrow
B_{R}\left(  o\right)  \backslash\left\{  o\right\}  $ is a Riemannian
isometry. To conclude, note that, by the assumptions on $g$, this isometry
smoothly extends even to the origin $0.$\medskip

\textbf{Step 4.} We now assume that $u$ is a solution of (\ref{local1}) on all
of $M$. In case $r_{-G}=+\infty$, then it follows directly from Step 3 that
$\Phi:M_{-G}^{m}\rightarrow M$ is a Riemannian isometry. Accordingly, in what
follows, we assume $r_{-G}<+\infty$.\medskip

\textbf{Step 5.} We show that $\partial B_{r_{-G}}\left(  o\right)  $ is
discrete, hence a finite set. Indeed, for every $x\in\partial B_{r_{-G}%
}\left(  o\right)  $, let $\gamma$ be a unit speed, minimizing geodesic from
$o$ to $x$. Then $\left\vert \nabla u\right\vert \circ\gamma\left(  t\right)
=Ag\left(  t\right)  \rightarrow0$ as $t\rightarrow r_{-G}$. Therefore,
$\partial B_{r_{-G}}\left(  o\right)  $ is made up by critical points of $u$.
Since $u$ satisfies the differential equation $\mathrm{Hess}\left(  u\right)
\left(  x\right)  =H\left(  r\left(  x\right)  \right)  u\left(  x\right)
\left\langle ,\right\rangle $ and, by assumption, $H\left(  r_{-G}\right)
\neq0$ and $u\neq0$ on $\partial B_{r_{-G}}\left(  o\right)  $, we deduce that
such critical points are non-degenerate (i.e., the quadratic form
$\mathrm{Hess}\left(  u\right)  $ has no zero eigenvalues) hence, by Morse
Lemma, they are isolated. Accordingly, $\partial B_{r_{-G}}=\left\{
p_{1},...,p_{k}\right\}  $, as claimed.\medskip

\textbf{Step 6. } We prove that $cut\left(  o\right)  =\left\{  O\right\}
=\partial B_{r_{-G}}\left(  o\right)  $, for some $O\in M$. Indeed, by Step 2,
the standard $m$-dimensional ball $\mathbf{B}_{r_{-G}}^{m}\left(  0\right)
\subset T_{o}M$ of radius $r_{-G}$ lies in the domain $D_{o}\subset T_{o}M$ of
the normal coordinates at $o.$ Therefore, it suffices to show that%
\begin{equation}
\exp_{o}\left(  \partial\mathbf{B}_{r_{-G}}^{m}\left(  0\right)  \right)
=\partial B_{r_{-G}}\left(  o\right)  =\left\{  O\right\}  . \label{global1}%
\end{equation}
If this occurs then $\partial\mathbf{B}_{r_{-G}}^{m}\left(  0\right)  $ is
precisely the tangential cut-locus of $o$ and, hence, $cut\left(  o\right)
=\left\{  O\right\}  $. Note that, in particular, all the geodesics issuing
from $o$ will meet at $O$ (and cannot minimizes distances past $r_{-G}$).

Now for the proof of (\ref{global1}). Let us observe that $\exp_{o}%
(\partial\mathbf{B}_{r_{-G}}^{m}\left(  0\right)  )\subseteq\overline
{B_{r_{-G}}\left(  o\right)  }$ and $\exp_{o}(\partial\mathbf{B}_{r_{-G}}%
^{m}\left(  0\right)  \cap D_{o})=\partial B_{r_{-G}}\left(  o\right)
\cap(M\backslash cut\left(  o\right)  )$. Since $B_{r_{-G}}\left(  o\right)  $
does not contain any cut-point of $o$, it follows that also the tangential cut
points in $\partial\mathbf{B}_{r_{-G}}^{m}\left(  0\right)  $ are mapped on
$\partial B_{r_{-G}}\left(  o\right)  $ by $\exp_{o}.$ Thus, $\exp
_{o}(\partial\mathbf{B}_{r_{-G}}^{m}\left(  0\right)  )=\partial B_{r_{-G}%
}\left(  o\right)  $. Now, recall from Step 5 that $\partial B_{r_{-G}}\left(
o\right)  $ is a finite set. Since $\partial\mathbf{B}_{r_{-G}}^{m}\left(
0\right)  $ is connected and $\exp_{o}$ is a continuous map, we conclude the
validity of (\ref{global1}).\medskip

\textbf{Step 7.} We note that $\Phi\left(  M_{-G}^{m}\right)  =M\backslash
\left\{  O\right\}  =B_{r_{-G}}\left(  o\right)  .$ Indeed, this follows
directly from Step 3 and Step 6.\medskip

\textbf{Step 8.} We finally deduce that $M$ is, diffeomorphically, a twisted
sphere. To this end, recall that, by Step 6, $M$ is compact. Moreover, $u$ is
a smooth function on $M$ with precisely two critical points, $o$ and $O$.
According to (\ref{local1}) and Step 5, these critical points are
non-degenerate. Therefore, to conclude, we can apply the (differentiable
version of) the classical result by G. Reeb.\medskip

This completes the proof of the Theorem.
\end{proof}

\section{Third order systems: from functions to (gradient) vector fields}

Recently, a lot of work has been made to characterize space-forms, and also
complex K\"{a}hler and quaternionic manifolds, via differential equations
involving vector fields instead of functions. We refer to \cite{EGKU-specific}
for a survey of such a kind of results. Let us focus the attention on
space-forms. It is a nice observation by Garcia-Rio, Kupeli and Unal,
\cite{GKU-spheres}, that if the vector field $Z$ on $M$ satisfies%
\begin{equation}
\left(  DDZ\right)  \left(  X,Y\right)  =k\left\langle Z,X\right\rangle Y,
\label{vector0}%
\end{equation}
for every vector fields $X,Y$ and for some constant $k\neq0,$ where $D$
denotes the covariant differentiation, then: (a) $Z$ has the special form%

\begin{equation}
Z=\frac{\nabla\operatorname{div}Z}{mk}, \label{vector2}%
\end{equation}
and, (b) the smooth function%
\[
u=\operatorname{div}Z
\]
satisfies%
\begin{equation}
\mathrm{Hess}\left(  u\right)  =ku\left\langle ,\right\rangle ,\text{ on }M.
\label{vector3}%
\end{equation}

Using this latter fact, the authors are able to reduce their characterizations
of space-forms to the \textquotedblleft scalar\textquotedblright\ cases
collected in Theorem \ref{th_spaceforms}. Note that, once we have chosen a
reference origin $o$ and used polar coordinates with respect to $o$, the
function $u$ turns out to be radial and hence, by (\ref{vector2}), $Z$ is a
radial gradient vector field.

One may therefore ask whether similar characterizations hold for a generic
model up to considering solutions of%
\begin{equation}
\left(  DDZ\right)  \left(  X,Y\right)  =K\left(  r\left(  x\right)  \right)
\left\langle Z,X\right\rangle Y, \label{vector1}%
\end{equation}
for a suitable smooth, real valued function $K\left(  t\right)  $, thus
extending Theorem \ref{th_global} to vector field equations. Inspection of
what happens on a generic model suggests that this is the case. Indeed,
suppose we are given a model $M_{-G}^{m}$ with corresponding warping function
$g$. In view of what we observed above, it is quite natural to consider the
radial, gradient vector field%
\[
Z_{x}=\nabla\left(  \int_{0}^{r\left(  x\right)  }y\left(  s\right)
ds+B\right)  =y\left(  r\left(  x\right)  \right)  \nabla r,
\]
where $B\in\mathbb{R}$ is an arbitrary constant. Straightforward calculations
show that%
\begin{align*}
\left\langle \left(  DDZ\right)  \left(  X,Y\right)  ,W\right\rangle  &
=y^{\prime\prime}dr\left(  X\right)  dr\left(  Y\right)  dr\left(  Z\right) \\
&  +y^{\prime}\mathrm{Hess}\left(  r\right)  \left(  X,Y\right)  dr\left(
W\right) \\
&  +y^{\prime}\mathrm{Hess}\left(  r\right)  \left(  X,W\right)  dr\left(
Y\right) \\
&  +y^{\prime}\mathrm{Hess}\left(  r\right)  \left(  Y,W\right)  dr\left(
X\right) \\
&  +y\left(  D_{X}\mathrm{Hess}\left(  r\right)  \right)  \left(  Y,W\right)
.
\end{align*}
On the other hand, using (\ref{radial2}), we see that%
\begin{align}
\left(  D_{X}\mathrm{Hess}\left(  r\right)  \right)  \left(  Y,W\right)   &
=\left\{  \frac{\left(  gg^{\prime}\right)  ^{\prime}}{gg^{\prime}}%
-2\frac{g^{\prime}}{g}\right\}  \mathrm{Hess}\left(  r\right)  \left(
Y,W\right)  dr\left(  X\right) \label{vector4}\\
&  -\frac{g^{\prime}}{g}\mathrm{Hess}\left(  r\right)  \left(  X,Y\right)
dr\left(  W\right)  -\frac{g^{\prime}}{g}\mathrm{Hess}\left(  r\right)
\left(  X,W\right)  dr\left(  Y\right)  ,\nonumber
\end{align}
holds for every vector fields $X,Y,W$. Whence, we deduce%

\begin{align}
\left\langle \left(  DDZ\right)  \left(  X,Y\right)  ,W\right\rangle  &
=y^{\prime\prime}dr\left(  X\right)  dr\left(  Y\right)  dr\left(  W\right)
\label{vector5}\\
&  +\left(  y^{\prime}-y\frac{g^{\prime}}{g}\right)  \mathrm{Hess}\left(
r\right)  \left(  X,Y\right)  dr\left(  W\right) \nonumber\\
&  +\left(  y^{\prime}-y\frac{g^{\prime}}{g}\right)  \mathrm{Hess}\left(
r\right)  \left(  X,W\right)  dr\left(  Y\right) \nonumber\\
&  +\left(  y^{\prime}+y\frac{\left(  gg^{\prime}\right)  ^{\prime}%
}{gg^{\prime}}-2y\frac{g^{\prime}}{g}\right)  \mathrm{Hess}\left(  r\right)
\left(  Y,W\right)  dr\left(  X\right)  .\nonumber
\end{align}
Since%
\begin{align}
K\left(  r\right)  \left\langle Z,X\right\rangle \left\langle
Y,W\right\rangle  &  =K\left(  r\right)  ydr\left(  X\right)  dr\left(
Y\right)  dr\left(  W\right) \label{vector6}\\
&  +K\left(  r\right)  y\frac{g}{g^{\prime}}\mathrm{Hess}\left(  r\right)
\left(  Y,W\right)  dr\left(  X\right)  ,\nonumber
\end{align}
it follows that (\ref{vector1}) is verified for the chosen vector fields $X,Y$
if and only if the following equation%
\begin{align}
0  &  =\left(  y^{\prime\prime}-Ky\right)  dr\left(  X\right)  dr\left(
Y\right)  dr\left(  W\right) \label{vector7}\\
&  +\left(  y^{\prime}-y\frac{g^{\prime}}{g}\right)  \mathrm{Hess}\left(
r\right)  \left(  X,Y\right)  dr\left(  W\right) \nonumber\\
&  +\left(  y^{\prime}-y\frac{g^{\prime}}{g}\right)  \mathrm{Hess}\left(
r\right)  \left(  X,W\right)  dr\left(  Y\right) \nonumber\\
&  +\left(  y^{\prime}+y\frac{\left(  gg^{\prime}\right)  ^{\prime}%
}{gg^{\prime}}-2y\frac{g^{\prime}}{g}-Ky\frac{g}{g^{\prime}}\right)
\mathrm{Hess}\left(  r\right)  \left(  Y,W\right)  dr\left(  X\right)
\nonumber
\end{align}
is satisfied for every $W$. Using appropriate choices of $X,Y,W$ we
immediately see that equation (\ref{vector1}) is equivalent to%
\begin{equation}
\left\{
\begin{array}
[c]{l}%
y^{\prime\prime}-Ky=0\\
y^{\prime}-yg^{\prime}/g=0\\
y^{\prime}+y\left(  gg^{\prime}\right)  ^{\prime}/gg^{\prime}-2yg^{\prime
}/g-Kyg/g^{\prime}=0.
\end{array}
\right.  \label{vector8}%
\end{equation}
Whence, up to imposing $y\left(  0\right)  =0$ (which is a natural assumption
in order to extend the above computations to the pole of $M_{-G}^{m}$) we
conclude that these conditions imply%
\begin{equation}
K\left(  r\right)  =G\left(  r\right)  \text{,\qquad}y\left(  r\right)
=Ag\left(  r\right)  , \label{vector9}%
\end{equation}
for any constant $A\neq0$. We have thus obtained the following

\begin{lemma}
\label{lemma_vector-model}A necessary and sufficient condition for equation
(\ref{vector1}) on $M_{-G}^{m}$ to possess a (non-trivial) radial, gradient
vector field solution $Z$ is that $K\left(  r\right)  =G\left(  r\right)  $.
In this case,
\begin{equation}
Z_{x}=\nabla\left(  A\int_{0}^{r\left(  x\right)  }g\left(  s\right)
ds+B\right)  , \label{vector10}%
\end{equation}
where $A\neq0$ and $B\in\mathbb{R}$ are arbitrary constants.
\end{lemma}

Observe that $Z$ is the gradient vector field associated to the radial
solution $u\left(  x\right)  =\alpha\left(  r\left(  x\right)  \right)  $ of
the \textquotedblleft scalar\textquotedblright\ equation (\ref{radial}). Also,
as we already remarked at the beginning of the section, if $G\left(  r\right)
\equiv k$ a non-zero constant then, according to (\ref{vector2}), any solution
$Z$ of (\ref{vector0}) must be of the form $Z=\nabla u$ where
$u=\operatorname{div}Z/mk,$ and equation (\ref{vector0}) becomes%
\[
\left(  D\mathrm{Hess}\left(  u\right)  \right)  \left(  X;Y,W\right)
=k\left\langle \nabla u,X\right\rangle \left\langle Y,W\right\rangle .
\]
According to these considerations, we are naturally led to state the next
rigidity result which represents a genuine extension of Theorem
\ref{th_spaceforms-vector} stated in the Introduction. Our approach is rather
different from that presented in \cite{GKU-spheres}. Indeed, the reduction
procedure outlined above cannot be carry over in this more general situation.

\begin{theorem}
\label{th_vector}Let $\left(  M,\left\langle ,\right\rangle \right)  $ be an
$m$-dimensional, complete Riemannian manifold, let $o\in M$ be a reference
origin and set $r\left(  x\right)  =\mathrm{dist}_{\left(  M,\left\langle
,\right\rangle \right)  }\left(  x,o\right)  $. A necessary and sufficient
condition for the existence of an isometric imbedding $\Phi:M_{-G}%
^{m}\rightarrow M$ is that there exists a non-trivial, smooth solution
$u:B_{r_{-G}}\left(  o\right)  \rightarrow\mathbb{R}$ \ of the problem%
\begin{equation}
\left\{
\begin{array}
[c]{l}%
\left(  D\mathrm{Hess}\left(  u\right)  \right)  \left(  X;Y,W\right)
=G\left(  r\left(  x\right)  \right)  \left\langle \nabla u,X\right\rangle
\left\langle Y,W\right\rangle \\
\mathrm{Hess}\left(  u\right)  \left(  o\right)  =A\left\langle ,\right\rangle
\\
\left\vert \nabla u\right\vert \left(  o\right)  =0,
\end{array}
\right.  \label{vector12}%
\end{equation}
for some $A\neq 0$.
Furthermore, if $u$ is a solution of (\ref{vector12}) on all of $M$, then the
following holds:

\begin{enumerate}
\item[(a)] If $r_{-G}=+\infty$, then $M$ is isometric to the model $M_{-G}%
^{m}$.

\item[(b)] In case $r_{-G}<+\infty$, and $g^{\prime}\left(  r_{-G}\right)
\neq0$, then $cut\left(  o\right)  =\left\{  O\right\}  $ for some $O\in M$,
and $\Phi\left(  M_{-G}^{m}\right)  =M\backslash\left\{  O\right\}  $.
Moreover, $M$ is diffeomorphically a twisted sphere.
\end{enumerate}
\end{theorem}

\begin{remark}
In case $G\left(  s\right)  \equiv k$ it can be shown that  assumption
$\mathrm{Hess}\left(  u\right)  \left(  o\right)  =A\left\langle
,\right\rangle $ is unessential. Furthermore, if $k<0$ then even the request
$\left\vert \nabla u\right\vert \left(  o\right)  =0$ \ can be omitted.
\end{remark}

Comparing Theorems \ref{th_global} and \ref{th_vector} we see that the
characterization of a model $M_{-G}^{m}$ \ via a third order system seems to
be more natural. Indeed, the system involves directly the radial sectional
curvature $-G\left(  r\right)  $ of the model. On the other
hand, in the situation of second order systems, we are able to characterize
the same space via a one-parameter family of differential systems as remarked
e.g. in Corollary \ref{cor_spaceforms}. This further characterizations are
invisible from the third order point of view.

\begin{proof}
Let us begin by showing that $u\left(  x\right)  =\alpha\left(  r\left(
x\right)  \right)  $ where%
\begin{equation}
\alpha\left(  t\right)  =A\int_{0}^{t}g\left(  s\right)  ds+B,\label{vector13}%
\end{equation}
for some constant  $B\in\mathbb{R}$. To this aim, let $x\in
B_{r_{-G}}\left(  o\right)  $ be fixed and let $\gamma\left(  s\right)
:[0,r\left(  x\right)  ]\rightarrow B_{r_{-G}}\left(  o\right)  $ be a unit
speed, minimizing geodesic from $\gamma\left(  0\right)  =o$ to $\gamma\left(
r\left(  x\right)  \right)  =x$. Then, evaluating (\ref{vector12}) along
$\gamma$, we readily deduce that%
\[
y\left(  s\right)  =u\circ\gamma\left(  s\right)
\]
solves the Cauchy problem%
\begin{equation}
\left\{
\begin{array}
[c]{l}%
y^{\prime\prime\prime}=G\left(  s\right)  y^{\prime}\\
y\left(  0\right)  =B\\
y^{\prime}\left(  0\right)  =0\\
y^{\prime\prime}\left(  0\right)  =A,
\end{array}
\right.  \label{vector14}%
\end{equation}
where %
$B=u\left(  o\right)$. 
Since $G=g^{\prime\prime}/g$, integrating (\ref{vector14}) we deduce that%
\[
y\left(  s\right)  =\alpha\left(  s\right)  .
\]
Evaluating this latter at $s=r\left(  x\right)  $ we conclude that $u\left(
x\right)  =\alpha\left(  r\left(  x\right)  \right)  $ as desired.

As in Step 2 of the proof of Theorem \ref{th_global}, it follows from Bishop
density result that%
\[
cut\left(  o\right)  \cap B_{r_{-G}}\left(  o\right)  =\emptyset.
\]
On the other hand, using equation (\ref{vector12}) we have that the Riemann
curvature tensor of $M$ satisfies%
\begin{align*}
\mathrm{Riem}\left(  W,X,\nabla u,Y\right)   &  =\left(  D\mathrm{Hess}\left(
u\right)  \right)  \left(  W;X,Y\right)  -\left(  D\mathrm{Hess}\left(
u\right)  \right)  \left(  X;W,Y\right) \\
&  =G\left(  r\left(  x\right)  \right)  \left\{  \left\langle \nabla
u,W\right\rangle \left\langle X,Y\right\rangle -\left\langle \nabla
u,X\right\rangle \left\langle W,Y\right\rangle \right\}  ,
\end{align*}
for every vector fields $X,Y,W$. Since $\nabla u=Ag\left(  r\right)  \nabla
r$, choosing $X=\nabla r$ and $W=Y\in\left(  \nabla r\right)  ^{\bot}$ we
deduce that the radial sectional curvature of $M$ is given by%
\begin{equation}
Sec_{rad}\left(  x\right)  =-G\left(  r\left(  x\right)  \right)  .
\label{vector15}%
\end{equation}
Therefore, by Hessian comparisons, \cite{PRS-geomanal},%
\begin{equation}
\mathrm{Hess}\left(  r\right)  =\frac{g^{\prime}}{g}\left\{  \left\langle
,\right\rangle -dr\otimes dr\right\}  ,\text{ on }B_{r_{-G}}\left(  o\right)
\backslash\left\{  o\right\}  . \label{vector16}%
\end{equation}
Now the proof can be easily completed following the arguments of Theorem
\ref{th_global}. Indeed, setting%
\[
\exp_{o}^{\ast}\left\langle ,\right\rangle =dr\otimes dr+\sigma_{\alpha\beta
}\left(  r,\theta\right)  \theta^{\alpha}\otimes\theta^{\beta},
\]
and using (\ref{vector16}) into (\ref{local5}) yields the validity of
(\ref{local6}) which, once integrated, gives%
\[
\sigma_{\alpha\beta}\left(  r,\theta\right)  =g\left(  r\right)  ^{2}%
\delta_{\alpha\beta}.
\]
We have thus established that $B_{r_{-G}}\left(  o\right)  $ is isometric to
$M_{-G}^{m}$. In particular, $u$ satisfies%
\begin{equation}
\mathrm{Hess}\left(  u\right)  \left(  x\right)  =Ag^{\prime}\left(  r\left(
x\right)  \right)  \left\langle ,\right\rangle \text{, on }B_{r_{-G}}\left(
o\right)  . \label{vector17}%
\end{equation}
Suppose now that $u$ is defined on all of $M$. In case $r_{-G}=+\infty$ we
immediately conclude that $M$ is isometric to $M_{-G}^{m},$ as stated in (a).
On the other hand, assume that $r_{-G}<+\infty$, hence $g\left(
r_{-G}\right)  =0,$ and $g^{\prime}\left(  r_{-G}\right)  \neq0$. Having fixed
$x\in\partial B_{r_{-G}}\left(  o\right)  $ and a unit vector $v\in T_{x}M$,
let $\gamma:[0,r_{-G}]\rightarrow M$ be a minimizing geodesic from
$\gamma\left(  0\right)  =o$ to $\gamma\left(  r_{-G}\right)  =x$. Obviously,
$\gamma\left(  t\right)  \in B_{r_{-G}}\left(  o\right)  $ for every
$t<r_{-G}$. Next, consider $v\left(  t\right)  $ the vector field obtained by
parallel transport of $v$ along $\gamma$. Then, according to (\ref{vector17}),%
\[
\mathrm{Hess}\left(  u\right)  \left(  x\right)  \left(  v,v\right)
=\lim_{t\rightarrow r_{-G}}\mathrm{Hess}\left(  u\right)  \left(
\gamma\left(  t\right)  \right)  \left(  v\left(  t\right)  ,v\left(
t\right)  \right)  =Ag^{\prime}\left(  r_{-G}\right)  \neq0,
\]
proving that $\partial B_{r_{-G}}$ is made up entirely by non-degenerate
critical points. Therefore, following exactly Steps 5--8 in the proof of
Theorem \ref{th_global} we conclude the validity of the global properties of
$M$ collected in (b).
\end{proof}

We conclude the section by stating the following

\begin{problem}
Suppose that the equation $\left(  DDZ\right)  \left(  X,Y\right)  =G\left(
r\left(  x\right)  \right)  \left\langle Z,X\right\rangle Y$, $G\neq0$, has a
non-trivial solution $Z$ with $Z_{o}=0$. Is $M$ isometric to $M_{-G}^{m}$? Is
it necessary to impose some further assumption on $M$?
\end{problem}

Observe that, even in this more general situation,%
\[
Z=\frac{\nabla\operatorname{div}Z}{mG\left(  r\right)  }.
\]
However, this time, it does not follow from this expression that $Z$ is
gradient. Needless to say, the reduction procedure of \cite{GKU-spheres}
cannot be applied directly in the present situation.


\bigskip

\begin{thebibliography}{9}                                                                                                %


\bibitem {Bishop-cutloci}R. L. Bishop, \textit{Decomposition of cut loci.}
Proc. Amer. Math. Soc. \textbf{65} (1977), 133--136.

\bibitem {EGKU-specific}F. Erkekoglu, E. Garcia-Rio, D. Kupeli, B. Unal,
\textit{Characterizing specific Riemannian manifolds by differential
equations. }Acta Appl. Math. \textbf{76} (2003), 195--219.

\bibitem {GKU-spheres}E. Garcia-Rio, D. Kupeli and B. Unal,\textit{ On a
differential equation characterizing Euclidean spheres. }J. Differential
Equations \textbf{194 }(2003), 287--299.

\bibitem {Kanai-structure}M. Kanai, \textit{On a differential equation
characterizing a Riemannian structure of a manifold.} Tokyo J. Math.
\textbf{6} (1983), 143--151.

\bibitem {Obata-sphere}M. Obata,\textit{ Certain conditions for a Riemannian
manifold to be iosometric with a sphere. }J. Math. Soc. Japan \textbf{14}
(1962), 333--340.

\bibitem {Petersen-riemannian}P. Petersen, \textit{Riemannian geometry.}
Graduate Texts in Mathematics, \textbf{171}. Springer-Verlag, New York, 1998.

\bibitem {PRS-geomanal}S. Pigola, M. Rigoli, A.G. Setti, \textit{Vanishing and
finiteness results in geometric analysis. A generalization of the Bochner
technique.} Progress in Mathematics, \textbf{266}. Birkh\"{a}user Verlag,
Basel, 2008.

\bibitem {Tashiro-vector}Y. Tashiro, \textit{Complete Riemannian manifolds and
some vector fields.} Trans. Amer. Math. Soc. \textbf{117} (1965), 251--275.

\bibitem {Wolter-cutloci}F.-E. Wolter, \textit{Distance function and cut loci
on a complete Riemannian manifold. }Arch. Math. (Basel) \textbf{32} (1979), 92--96.
\end{thebibliography}
\end{document}